# *Existence and uniqueness theorem for ODE: an overview*

*Swarup Poria and Aman Dhiman*
*Department of Applied Mathematics,*
*University of Calcutta,*
*92, A.P.C.Road, Kolkata-700009, India*

**Abstract:** The study of existence and uniqueness of solutions became important due to the lack of general formula for solving nonlinear ordinary differential equations (ODEs). Compact form of existence and uniqueness theory appeared nearly 200 years after the development of the theory of differential equation. In the article, we shall discuss briefly the differences between linear and nonlinear first order ODE in context of existence and uniqueness of solutions. Special emphasis is given on the Lipschitz continuous functions in the discussion.

**1. Introduction:** Differential equations are essential for a mathematical description of nature, many of the general laws of nature-in physics, chemistry, biology, economics and engineering –find their most natural expression in the language of differential equation. Differential Equation(DE) allows us to study all kinds of evolutionary processes with the properties of determinacy; finite-dimensionality and differentiability. The study of DE began very soon after the invention of differential and integral calculus. In 1671, Newton had laid the foundation stone for the study of differential equations. He was followed by Leibnitz who coined the name differential equation in 1676 to denote relationship between differentials $dx$ and $dy$ of two variables $x$ and $y$. The fundamental law of motion in mechanics, known as Newton's second law is a differential equation to describe the state of a system. Motion of a particle of mass $m$ moving along a straight line under the influence of a specified external force $F(t, x, x')$ is described by the following DE

$$mx'' = F(t, x, x'); \quad \left(x' = \frac{dx}{dt}, x'' = \frac{d^2x}{dt^2}\right) \qquad (1)$$

At early stage, mathematicians were mostly engaged in formulating differential equations and solving them but they did not worry about the existence and uniqueness of solutions.

More precisely, an equation that involves derivatives of one or more unknown dependent variables with respect to one or more independent variables is known as differential equations. An equation involving ordinary derivatives of one or more dependent variables with respect to single independent variable is called an ordinary differential equations *(ODEs).* A general form of an ODE containing one independent and one dependent variable is $F(x, y, y', y'', \ldots\ldots, y^n) = 0$ where F is an arbitrary function of $x, y, y', \ldots\ldots, y^n$, here $x$ is the independent variable while y being the dependent variable and $y^n \equiv \frac{d^n y}{dx^n}$. The order of an ODE is the order of the highest derivative appearing in it. Equation (1) is an example of second order ODE. On the other hand, partial differential equations are those which have two or more independent variables.

Differential equations are broadly classified into linear and non-linear type. A differential equation written in the form $F(x, y, y', y'', \ldots, y^n) = 0$ is said to be linear if $F$ is a linear function of the variables $y, y', y'', \ldots, y^n$ (except the independent variable $x$). Let $L$ be a differential operator defined as

$$L \equiv a_n(x)D^n + a_{(n-1)}(x)D^{n-1} + \cdots + a_1(x)D + a_0(x)$$

where $D^n \equiv \frac{d^n y}{dx^n}$; (for n=1,2,....). One can easily check that the operator $L$ satisfies the condition for linearity i.e.

$$L(ay_1 + by_2) = aL(y_1) + bL(y_2)$$

for all $y_1, y_2$ with scalars $a$ and $b$. Therefore a linear ODE of order $n$ can be written as

$$a_n(x)y^{(n)} + a_{n-1}(x)y^{(n-1)} + \cdots + a_0(x)y = f(x) \quad (2)$$

An ODE is called nonlinear if it is not linear. In absence of damping and external driving, the motion of a simple pendulum (Fig.1) is governed by

$$\frac{d^2\theta}{dt^2} + \frac{g}{L}\sin\theta = 0 \quad (3)$$

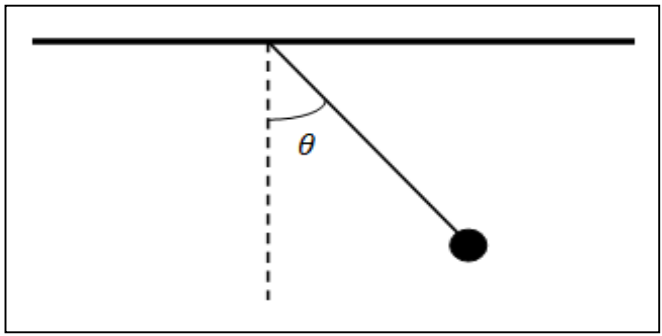

**Figure 1:** Oscillating simple pendulum

where $\theta$ is the angle from the downward vertical, $g$ is the acceleration due to gravity and $L$ is the length of the pendulum. This is a famous example of nonlinear ODE.

In this article, we shall confine our discussion to single scalar first order ODE only. An initial value problem(IVP) for a first order ODE is the problem of finding solution $y = y(x)$ that satisfies the initial condition $y(x_0) = y_0$ where $x_0, y_0$ are some fixed values. We write the IVP of first order ODE as

$$\frac{dy}{dx} = f(x, y), \quad y(x_0) = y_0. \quad (4)$$

In other words, we intend to find a continously differentiable function $y(x)$ defined on a confined interval of $x$ such that $\frac{dy}{dx} = f(x, y(x))$ and $y(x_0) = y_0$. Such a function $y = y(x)$ is called a solution to the equation (4). In integral form $y = y(x)$ is the solution of equation (4) if it satisfies the integral equation $y(x) = y_0 + \int_{x_0}^{x} f(u, y(u))du$.

There are several well known methods namely separation of variables, variation of parameters, methods using integrating factors etc for solving first order ODEs. In this section, we discuss briefly a method for solving first order linear differential equations. Any first order linear equation can be written as

$$\frac{dy}{dx} + p(x)y = q(x). \quad (5)$$

We assume that $p(x)$ and $q(x)$ are continuous functions of $x$ on the interval $a < x < b$. We look for a function $\mu(x) \neq 0$ such that

$$\mu(x)[y' + p(x)y] = \frac{d}{dx}[\mu(x)y] = \mu(x)y' + \mu'(x)y.$$

Then $\mu(x)$ must satisfy $\mu(x)p(x)y = \mu'(x)y$ and therefore $\mu(x) = \exp[\int p(x)\,dx]$. The general solution of equation (5) is given as
$$y\mu(x) = \int q(x)\mu(x)dx + c,$$
where c is an arbitrary constant. Therefore, it is obvious that all linear first order ODEs with continous coefficients $p(x)$ and $q(x)$ can be solved exactly.

In this context it is important to discuss about the following linear ODE:
$$y' + \frac{2}{x}y = 4x \quad \text{on } 0 < x < \infty. \tag{6}$$
Notice that $y = x^2 + \frac{c}{x^2}$ (*where c is arbitrary constant*) is a solution of equation (6). Setting the initial condition as $y(1) = 2$ the solution becomes $y = x^2 + \frac{1}{x^2}$. It is important to note that the solution becomes unbounded as $x \to 0$. This is not surprising since $x = 0$ is a point of dicontinuity of coefficient of $y$. However, if we choose initial condition as $y(1) = 1$ the solution is $y = x^2$ is well behaved as $x \to 0$. Therefore linear ODE with discontinuous coefficients are not easy and deeper understanding is required to handle these kinds of equations.

**2. Existance and uniqueness theorem:** From the above discussions it is quite clear that first order linear ODEs are easy to solve because there are straight forward methods for finding solutions of such equations *(with continuous coefficients)*. In contrast, there is no general method for solving even a first order nonlinear ODEs. Some typical nonlinear first order ODEs can be solved by Bernoulli's method(1697), method of seperation of variables, method of variation of parameters, method using integrating factors etc. The lack of general formula for solving nonlinear ODEs has two important consequences. Firstly, methods which yield approximate solution *(numerical)* and give qualitative information about solutions assume greater significance for nonlinear equations. Secondly, questions dealing with the existence and uniqueness of solutions became important.

In this article, we shall discuss briefly about existence and uniqueness of solutions of first oder ODEs. In ODE theory the following questions are naturally arising :

- Given an IVP is there a soultion to it *(question of existence)*?
- If there is a solution is the solution unique *(question of uniqueness)*?
- For which values of x does the solution to IVP exists *(the interval of existence)*?

The fundamentally important question of existence and uniqueness of solution for IVP was first answered by Rudolf Lipschitz in 1876 *(nearly 200 years later than the development of ODE)*. In 1886 Giuseppe Peano discovered that the IVP (4) has a solution *(it may not be unique)* if $f$ is a continous function of $(x, y)$. In 1890 Peano extended this theorem for system of first order ODE using method of successive approximation. In 1890 Charles Emile Picard and Ernst Leonard Lindelöf presented existence and uniqueness theorem for the solutions of IVP (4). According to Picard-Lindelöf theorem if $f$ and $\frac{\partial f}{\partial y}$ are continous functions of x, y in some rectangle: $\{(x, y): \alpha < x < \beta; \gamma < y < \delta\}$ containing the point $(x_0, y_0)$ then in some interval $x_0 - \delta < x < x_0 + \delta$ ($\delta > 0$) there exists a unique solution of IVP (4). In 1898

W.F.Osgood discovered a proof of Peano's theorem in which the solution of the differential equation is presented as the limit of a sequence of functions(*using Arzela-Ascoli theorem*).

Let us state the existence and uniqueness theorem which Lipschitz knew. Let $f(x,y)$ in (4) be continuous on a rectangle $D = \{(x,y): x_0 - \delta < x < x_0 + \delta; y_0 - b < y < y_0 + b\}$ then there exists a solution in D, furthermore, if $f(x,y)$ is Lipschitz continuous with respect to $y$ on a rectangle $R$ (*possibly smaller than D*) given as $R = \{(x,y): x_0 - a < x < x_0 + a; y_0 - b < y < y_0 + b, a < \delta\}$ (Fig.2) the solution in $R$ shall be unique.

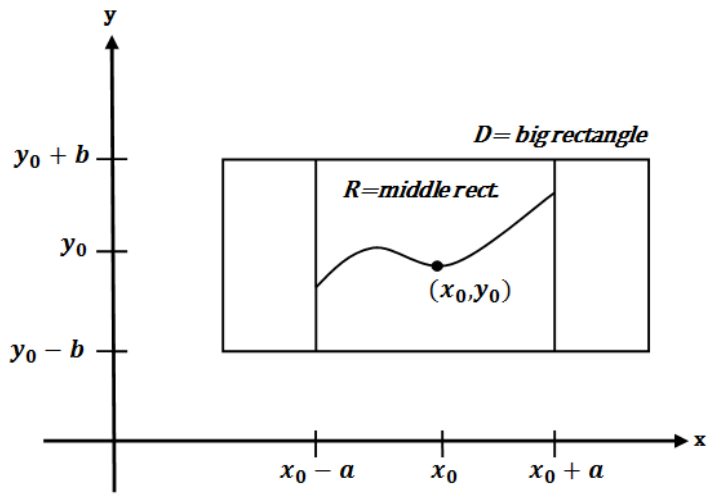

**Figure 2:** *Showing initial condition ($x_0, y_0$) enclosed in neighbourhoods R and D*

There are many ways to prove existence and uniqueness theorem one elegant proof of this theorem is given by using Banach fixed point theorem[1]. The proof is beyond the scope of the article. Lipschitz continuity was used by Lipschitz to prove the existence and uniqueness of solutions to IVP of ODE in 1876.

**2.1. Lipschitz continuity (local and global):** Understanding Lipschitz continuity is necessary to realize existence and uniqueness theory of ODE. A function $f(x,y)$ is said to be locally Lipschitz or locally Lipschitz continuous at a point $y_0 \in D$ (*open and connected set*) if $y_0$ has a neighbourhood $D_0$ such that $f(x,y)$ satisfies

$$|f(x,y_1) - f(x,y_2)| < L|y_1 - y_2| \quad \forall\, (x,y_1), (x,y_2) \in D_0$$

where L being fixed over the neighborhood $D_0$. A function is said to be locally Lipschitz in a domain D if it is locally Lipschitz at every point of the domain D. We denote the set of all locally Lipschitz functions by $\mathcal{L}_l$ and symbolically we say $f \in \mathcal{L}_l$ if $f$ is locally Lipschitz. Moreover, a function is said to be globally Lipschitz (*or $f \in$ set of globally Lipschitz function, $\mathcal{L}_g$*) if

$$|f(x,y_1) - f(x,y_2)| < L|y_1 - y_2| \quad \forall\, (x,y_1), (x,y_2) \in D$$

with the same Lipschitz constant L for entire $D$. Geometrically, the slope of any chord joining two points in the neighborhood of $y_0$ is bounded then we say that $f$ is locally Lipschitz continuous at $y_0$, if this slope is bounded in the neighborhood of every point of domain D we say that $f$ is locally Lipschitz in D. Furthermore, if this bound remains fixed over the entire domain then the function is globally Lipschitz and $f \in \mathcal{L}_g$.

In the proceeding discussion, whenever the function $f(x,y)$ is a function of $y$ only we write it as $f(y)$.

The difference between $\mathcal{L}_l$ and $\mathcal{L}_g$ functions can be understood from the following discussion. Let us consider the function $f(y) = y^2$ on $\mathbb{R}$. For $y_0 \in \mathbb{R}$ we observe that

$$\sup_{y \in (y_0-1, y_0+1)} |f'(y)| = \sup_{y \in (y_0-1, y_0+1)} |2y| \leq 2|y_0| + 1 \qquad (7)$$

as $|y| \leq |y_0| + 1$ for $y \in (y_0 - 1, y_0 + 1)$ by triangle inequality. Now picking up two points $y_1, y_2 \in (y_0 - 1, y_0 + 1)$ it follows by mean value theorem that for some $\xi \in (y_1, y_2)$, we have

$$|f(y_1) - f(y_2)| = |f'(\xi)(y_1 - y_2)|$$
$$\leq \sup_{\theta \in (y_0-1, y_0+1)} |f'(\theta)||y_1 - y_2|$$

Using the result obtained in (7) we get

$$|f(y_1) - f(y_2)| \leq (2|y_0| + 1)|y_1 - y_2|$$

for all $y_1, y_2 \in (y_0 - 1, y_0 + 1)$. Hence the Lipschitz constant $L = (2|y_0| + 1)$ for $f$ on $(y_0 - 1, y_0 + 1)$, so $f \in \mathcal{L}_l$ for all $y_0$ of $\mathbb{R}$. Here we should observe that the Lipschitz constant L depends on $y_0$. In particular, if $y_0 \to \infty$ then $L \to \infty$. This is an indication that the function is not globally Lipschitz on $\mathbb{R}$ although it is locally Lipschitz on any interval. Indeed, for any $y_0 \neq 0$

$$\left|\frac{f(y_0) - f(0)}{y_0 - 0}\right| = |y_0| \to \infty \text{ as } y_0 \to \infty$$

This means there is no L that satisfies the global Lipschitz property for $y_0 \in \mathbb{R}$, hence $f \notin \mathcal{L}_g$.

Existence and uniqueness theory for linear IVP (5) is simple and straight forward. To check this, equation (5) is rewritten as

$$y' = -p(x)y - q(x) = f(x, y).$$

Now, $|f(x, y_1) - f(x, y_2)| = |-p(x)(y_1 - y_2)|$
$$< M|y_1 - y_2|$$

where $M = \max_{a \leq x \leq b} |p(x)|$. Here continuity of $p(x)$ for $x \in [a, b]$ is enough for existence of a bound M for all $x \in [a, b]$ and so $f \in \mathcal{L}_g$. In fact any function $f(x, y)$ of the form $a(x)y + b(x)$ is Lipschitz continuous in any closed interval. This shows that IVP of first order linear systems where $p(x)$ and $q(x)$ are continuous functions of $x$ has a unique solution in any interval [a,b] but this may not always be the case for nonlinear ODEs.

Let us start our discussion with a simplest nonlinear function namely piecewise linear function with a jump discontinuity. Any jump discontinuity makes the slope of the chord joining the two neighboring points of the point of discontinuity (one left point and one right point) unbounded and so the function is non-Lipschitz at the point of discontinuity. Let us consider the function

$$f(y) = \begin{cases} y & ; \quad y < a \\ y + k & ; \quad y \geq a \end{cases}$$

for $k \neq 0$ having point of discontinuity at $y = a$. Let A be a point on the left side of the point of discontinuity and C be a point on the right side of discontinuity point (Fig.3).

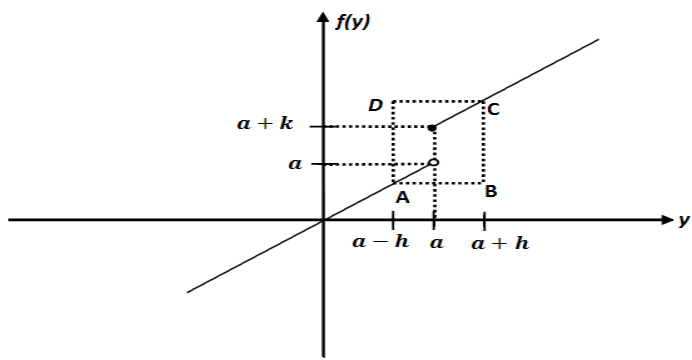

**Figure 3:** Shows graph of **f(y)** with discontinuity at **y = a**. Point A being left of discontinuity while B right of discontinuity.

Then the slope of $AC = \frac{BC}{AB} = \frac{(a+h+k)-(a-h)}{a+h-(a-h)} = \frac{2h+k}{2h} = 1 + \frac{k}{2h}$

Now as $h \to 0$ slope of $AC \to \infty$. Showing $f \notin \mathcal{L}_l$ in any neighborhood of $y = a$.

An obvious question arises at this juncture that does continuity ensure Lipschitz and vice-versa? Answer to the first part is NO!

We can easily see this after having a close inspection of the function $f(y) = \sqrt[3]{y}, y \in [-3,3]$. The function is continuous on [-3, 3] (see Fig.4) but in any small neighborhood of y=0 the slope of chord is unbounded as

$$\left|\frac{f(y_1)-f(y_2)}{y_1-y_2}\right| = \left|\frac{(y_1^{1/3}-y_2^{1/3})}{y_1-y_2}\right| = \left|\frac{1}{y_1^{2/3}+y_2^{2/3}+y_1^{2/3}y_2^{2/3}}\right| \to \infty \text{ as } y_1, y_2 \to 0.$$

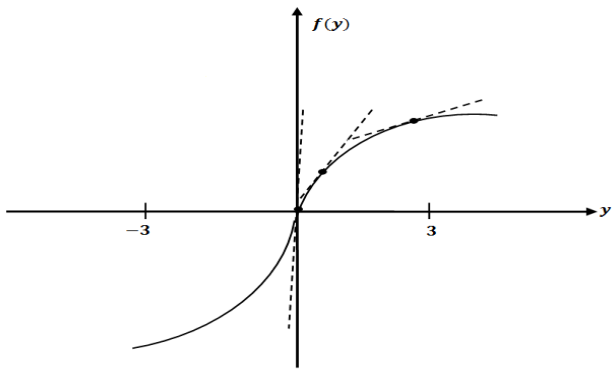

**Figure 4:** Graph of $f(y) = \sqrt[3]{y}$ showing continuity at origin

Therefore, all continuous functions are not necessarily Lipschitz continuous. Conversely, Lipschitz condition actually guarantees uniform continuity which is a stronger condition than just continuity. This easily follows from the definition of $\mathcal{L}_l$ functions. To show this let us consider a function $f \in \mathcal{L}_l$ on an interval $I$ then by defination

$$|f(y_1) - f(y_2)| < L|y_1 - y_2| \quad \forall \, y_1, y_2 \in I \quad (8)$$

Now $\forall \, \varepsilon > 0$ defining $\delta = \varepsilon/L$ equation (8) can be rewritten in a more familiar manner as

$$|f(y_1) - f(y_2)| < \varepsilon, \quad \forall \, |y_1 - y_2| < \delta$$

$y_1, y_2 \in I$. This is the condition required for uniform continuity. However, it is not definite that a uniformly continuous function is Lipschitz continuous. This can be observed in study of function $f(y) = \sqrt[3]{y}, y \in [-3,3]$. Here $f$ is continuous on the compact set $[-3,3]$, hence it is uniformily continuous on $[-3,3]$ although $f$ is not Lipschitz continuous there.

At this instant the reader might search for a relation between differentiable functions and Lipschitz continuous functions. Does Lipschitz continuity ensure differentiability and vice versa? Answer to both these questions is "NO". To justify this point, it is worth noticing that for the function $f(y) = |y|$ the derivative at $y = 0$ does not exist, however, the slope of the chord remains bounded about it 0. Since,

$$|f(y_1) - f(y_2)| = ||y_1| - |y_2|| \leq |y_1 - y_2| \quad \forall \, y_1, y_2 \in \mathbb{R}$$

Hence, $\frac{|f(y_1)-f(y_2)|}{|y_1-y_2|} < 1 = L$.

So $f \in \mathcal{L}_g$ with Lipschitz constant $L = 1$ throughout $\mathbb{R}$ although $f$ is not differentiable at $y = 0$. On the other hand, it is amazing that differentiability of $f$ is still not enough for Lipschitz condition to hold. For example let us consider the function

$$f(y) = \begin{cases} y^2 \sin\frac{1}{y^2}, & y \neq 0 \\ 0, & y = 0 \end{cases}.$$

This function is differentiable everywhere and its derivative is given as

$$f'(y) = \begin{cases} 2y \sin\frac{1}{y^2} - \frac{2}{y}\cos\frac{1}{y^2}, & y \neq 0 \\ 0, & y = 0 \end{cases}.$$

However, $|f'(y)| \to \infty$ as $y \to 0$, so the derivative is not continuous at $y = 0$. Now, to check for Lipschitz condition, we consider two sequence of points $\{y_{1,n}\} = \{(2n\pi + \frac{\pi}{2})^{-1/2}\}$ and $\{y_{2,n}\} = \{(2n\pi)^{-1/2}\}$ $(n = 1,2,\dots)$ which have the limit $y = 0$.

Now

$$\left|\frac{f(y_{1,n}) - f(y_{2,n})}{y_{1,n} - y_{2,n}}\right| = \frac{(2n\pi + \frac{\pi}{2})^{-1}}{(2n\pi)^{-1/2} - (2n\pi + \frac{\pi}{2})^{-1/2}}$$

$$= 4n\left((2n\pi)^{-1/2} + (2n\pi + \pi/2)^{-1/2}\right) \to \infty$$

as $n \to \infty$. Hence, $f \notin \mathcal{L}_l$.

With this example we realize that being just differentiable isn't sufficient for a function to be Lipschitz. However, when we impose a condition that the derivative is bounded or simply $f$ is continuously differentiable $(C^1)$ on a closed interval then $f$ turns out to be Lipschitz [4]. To show this let us suppose that $f$ is a continuously differentiable function $(f \in C^1)$ on $[a, b]$. This implies $f'$ is continuous on $[a, b]$ and so must be bounded there i.e. $|f'(y)| < M$, $y \in [a, b]$. Now for $y_1, y_2 \in [a, b]$ we have

$$|f(y_1) - f(y_2)| = \left|\int_{y_1}^{y_2} f'(u) du\right| \leq \int_{y_1}^{y_2} |f'(u)| \, du < M|y_1 - y_2|.$$

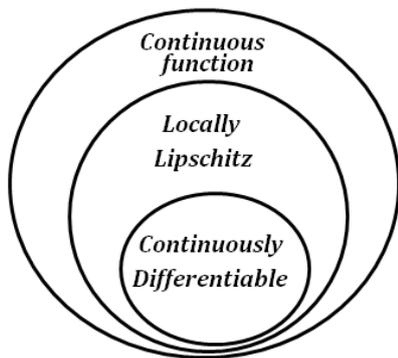

Therefore, $f$ is Lipschitz with Lipschitz constant $L = M$. Summarizing the results obtained from previous examples we have - every locally Lipschitz function is uniformly continuous and every differentiable function is not Lipschitz continuous but every continuously differentiable function is locally Lipschitz . A pictorial conclusion is represented in figure alonside.

Let us now see what happens for IVP (4) when $f(x,y)$ is non-Lipschitz but continuous. For this we consider

$$\frac{dy}{dx} = y^{1/3} \; ; \quad y(0) = 0. \tag{9}$$

Here $f(x,y) = y^{1/3}$ is non Lipschitz at zero *(as seen earlier)*. Solving it subject to the given initial condition we arrive at the solution $y(x) = \left(\frac{2x}{3}\right)^{3/2}$. Also it is easy to see $y(x) \equiv 0$ is also a solution of (9). Using the translation invariance w.r.t $x$ of the equation (9), we can patch together the solutions to obtain an infinite family of solutions. So for any $c \geq 0$

$$y(x) = \begin{cases} \left(\frac{2(x-c)}{3}\right)^{3/2}, & x \geq c \\ 0, & x < c \end{cases}$$

serves as a solution to this initial value problem. Hence there are infinitely many solutions of this IVP for infinitely many choices for c. This can be easily observed from the figure 5, here we have plotted $y(x)$ for $c_1 = 0, c_2 = 1, c_3 = 2$. Moreover, it is important to note that this IVP has uncountable infinite number of solutions [2].

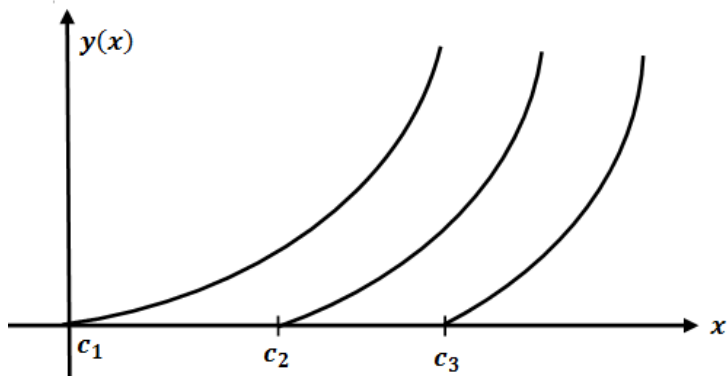

**Figure 5:** *Graph representing three solutions for IVP (9) with $c_1 = 0, c_2 = 1, c_3 = 2$*

It is worth mentioning here that Peano's existence condition as well as Lipschitz condition are respectively the sufficient conditions for the existence and uniqueness of solutions for IVP (4) but not necessary. This can be easily understood from the following IVP :

$$\frac{dy}{dx} = \frac{1}{y^2}, \quad y(x_0) = 0. \tag{10}$$

Here, $f(x,y) = \frac{1}{y^2}$ is discontinuous at $y = 0$ i.e. it fails to satisfy Peano's existence condition. However, a unique solution given as $y(x) = \left(3(x-x_0)\right)^{(1/3)}$ exists for every initial point chosen on the $x$-axis.

It is interesting to see that with the help of existence and uniqueness theorem we can predict the interval of existence of solution also. Considering the example

$$\frac{dy}{dx} = y^2, \quad y(0) = y_0.$$

Solution of this IVP is $y(x) = \frac{1}{\frac{1}{y_0} - x}$ which is defined on the interval $x \in [0, \frac{1}{y_0})$. Therefore the interval of existence for the nonlinear IVP is $[0, \frac{1}{y_0})$.

In fact if the function $f(x,y)$ is globally Lipchitz with respect to $y$ then solution remains defined globally. Like the case of $\frac{dy}{dx} = sinx$ with any initial condition solution stays uniquely defined for all $x \in \mathbb{R}$. Another important point to note that if $f(x,y)$ in (4) fails to be Lipschitz with respect to $y$ but is Lipschitz with respect to $x$ and $f(x_0, y_0) \neq 0$ then the IVP (4) has a unique solution if and only if the IVP $\frac{dx}{dy} = \frac{1}{f(x,y)}, y(x_0) = y_0$ has a unique solution [3]. Let us discuss it with the help of the IVP $\frac{dy}{dx} = \cos x + x\sqrt[3]{y}$ with $y(0) = 0$. Here $f(x,y)$ does not satisfy Lipschitz condition with respect to $y$. Since,

$$\left|\frac{f(x,y_1)-f(x,y_2)}{y_1-y_2}\right| = \left|\frac{x\left(y_1^{1/3}-y_2^{1/3}\right)}{y_1-y_2}\right|$$

$$= \left|\frac{x}{y_1^{2/3}+y_2^{2/3}+y_1^{2/3}y_2^{2/3}}\right| \to \infty \text{ as } y_1, y_2 \to 0$$

showing the slope of chord joining points $y_1, y_2$ in neighborhood of (0,0) is unbounded hence $f \notin \mathcal{L}_l$ with respect to $y$. While $\frac{\partial f}{\partial x} = -sinx + \sqrt[3]{y}$ is bounded in neighborhood of (0, 0) hence $f \in \mathcal{L}_l$ with respect to $x$. Therefore Lipschitz continuity with respect to $x$ may also be useful to establish uniqueness of solution.

**2.2. Algebraic properties of Lipschitz continuous functions:** We have so far realized the importance of Lipschitz continuity in the existence and uniqueness theory of IVP of ODE. We now discuss about some algebraic properties of $\mathcal{L}_l$ functions. Suppose $f_1, f_2 \in \mathcal{L}_l$ with Lipschitz constants $L_1, L_2$ on an interval I. Then for two points $y_1, y_2 \in I$ we have

$|(f_1+f_2)(y_1) - (f_1+f_2)(y_2)|$

$$= |f_1(y_1)+f_2(y_1)-f_1(y_2)-f_2(y_2)|$$
$$\leq |f_1(y_1)-f_1(y_2)| + |f_2(y_1)-f_2(y_2)|$$
$$\leq L_1|y_1-y_2| + L_2|y_1-y_2|$$
$$= (L_1+L_2)|y_1-y_2|$$

Hence $f_1 + f_2 \in \mathcal{L}_l$ with Lipschitz constant as $L_1 + L_2$. The argument is also valid for the difference of two functions. The Lipschitz constant here is again $L_1 + L_2$. Since,

$|(f_1-f_2)(y_1) - (f_1-f_2)(y_2)|$

$$= |f_1(y_1)-f_2(y_1)-f_1(y_2)+f_2(y_2)|$$
$$\leq |f_1(y_1)-f_1(y_2)| + |f_2(y_2)-f_2(y_1)|$$
$$\leq L_1|y_1-y_2| + L_2|y_2-y_1|$$
$$= (L_1+L_2)|y_1-y_2|$$

Hence $f_1 - f_2 \in \mathcal{L}_l$.

Furthermore, it is easy to see that $f \in \mathcal{L}_l$ implies $cf \in \mathcal{L}_l$ (*where c is a constant*) with Lipschitz constant as $|c|L$. Making use of the previous two properties we arrive at a general result. Suppose $f_1, f_2, \ldots f_n \in \mathcal{L}_l$ with Lipschitz constants $L_1, L_2 \ldots L_n$ on interval I then their linear combination, $c_1f_1 + c_2f_2 + \ldots\ldots\ldots + c_nf_n \in \mathcal{L}_l$ with the Lipschitz constant $|c_1|L_1 + |c_2|L_2 + \ldots\ldots + |c_n|L_n$.

Moreover, if $f_1, f_2 \in \mathcal{L}_l$ on an interval I then their product $f_1 f_2 \in \mathcal{L}_l$ on the interval I. $f_1/f_2 \in \mathcal{L}_l$ if $f_2$ is bounded below by a positive constant on I. In addition, if $f_1 \in \mathcal{L}_l$ on interval $I_1$ with Lipschitz constant $L_1$ and $f_2 \in \mathcal{L}_l$ on interval $I_2$ with Lipschitz constant $L_2$ and $f_1(I_1) \subset I_2$ then their composition $f_2 \circ f_1 \in \mathcal{L}_l$ on $I_1$ with Lipschitz constant $L_1L_2$. We can see this by considering $y_1, y_2 \in I_1$ we have

$$|f_2(f_1(y_1)) - f_2(f_1(y_2))| \leq L_2|f_1(y_1) - f_2(y_2)|$$
$$\leq L_1L_2|y_1-y_2|$$

showing $f_2 \circ f_1 \in \mathcal{L}_l$ with Lipschitz constant is $L_1L_2$. These algebraic properties are useful in analyzing existence and uniqueness of solutions for complicated ODEs.

**3. Conclusion:** This review article mainly focused to understand the Lipschitz condition and its connection with existence and uniqueness of solutions of IVP for ODE. We had seen through examples that Lipschitz condition guaruntees uniform continuity but it does not ensure differentiability of the function. It was shown that, continuity is sufficient for existence of solution and locally Lipschitz is a sufficient condition for uniqueness of the solution of a IVP of first order ODE. It was pointed out that a unique solution exists for IVP of a linear ODE in any interval where coefficients are continuous and bounded. Moreover we can solve it explicitly. On the other hand, issues are not so straight forward in case of nonlinear ODEs.

In the article our attention was primarily around a single IVP of a first order ODE. The generalization of result for a system of $n$–first order IVP is straight forward. Say, if we have a system of $n$- first order ODEs as shown in (11) :

$$\left.\begin{aligned} \frac{dy_1}{dx} &= f_1(x, y_1, y_2, \ldots, y_n) \\ \frac{dy_2}{dx} &= f_2(x, y_1, y_2, \ldots, y_n) \\ &-----\\ &-----\\ &-----\\ \frac{dy_n}{dx} &= f_n(x, y_1, y_2, \ldots, y_n) \end{aligned}\right\} \quad (11)$$

with initial conditions $y_1(x_0) = y_{10}, \; y_2(x_0) = y_{20}\ldots, y_n(x_0) = y_{n0}$. Here to examine the existence and uniqueness of solution for the system (11) we need to check the applicability of Lipschitz condition on function $f_i$ with respect to the variable $y_i$ ($i = 1, 2, \ldots, n$) in neighborhood about $y_{i0}$. This result can be further applied to check the existence and uniqueness of solution to a $n$-order equation by transforming it to a system of $n$-first order ODEs. This can be achieved by defining variables as $y' = y_1, y'_1 = y_2, \ldots \ldots \ldots, y'_{n-1} = y_n$. The differences between linear and nonlinear ODE in the context of existence and uniqueness were discussed here briefly. Interested readers may read nonlinear dynamics books [5] to see many more interesting features of nonlinear ODEs. The study of nonlinear DEs is still a very active area of research due to its application in modeling various physical, chemical, biological, engineering and social systems.